\documentclass{article}
\usepackage{latexsym,amsthm,amsfonts,epsfig,psfrag }
\begin{document}

\newtheorem{cor}{Corollary}
\newtheorem{lemma}{Lemma}
\newtheorem{theorem}{Theorem}
\newtheorem{prop}{Proposition}
\renewcommand{\proofname}{Proof}

\def \H{{I\!\!H}}
\def \E{{I\!\!E}}
\def \Z{{\mathbb Z}}

\begin{center}
{\LARGE  \bf
Spherical simplices  \\
\vspace{2pt}
generating discrete reflection groups}

\vspace{9pt}
{\large A.~Felikson\footnote{Supported
by Schweizerischer Nationalfonds No. 20-61379.00.}}
\end{center}

\vspace{10pt}

\begin{center}
\parbox{10cm}
{\small
{\bf Abstract.}
Let $P$ be a simplex in  $S^n$ and
 $G_P$ be a group generated by
the reflections with respect to the facets of $P$.
We are interested in the case when  the group $G_P$
is discrete.
In this case we say that $G$ generates the discrete reflection group  $G_P$.
We develop the criteria for a simplex generating a discrete 
reflection group.
We also describe all indecomposable spherical 
simplices generating discrete 
reflection groups.

}

\end{center}

\vspace{8pt}

\section*{Introduction}

\noindent
{\bf 1. }
Let $P$ be a convex polyhedron in the spherical space $S^n$,
Euclidean space $\E^n$ or hyperbolic space $\H^n$.
 Consider a group $G_P$ generated by
the  reflections with respect to the facets of $P$.
We call $G_P$  a {\bf reflection group generated by $P$}.
The problem is to list the polyhedra generating discrete reflection groups.

Although the problem is old and easy to state, the answer is known only for
some very simple combinatorial types of polyhedra.
Already in 1873, Schwarz \cite{Shwarz}  listed the spherical triangles
generating discrete groups.
In  1998, E.~Klimenko and M.~Sakuma \cite{Pink}
solved the problem
for hyperbolic triangles.
In \cite{polygons}, \cite{pyrprism}, \cite{simplex},
\cite{lambert} 
the problem was solved for hyperbolic quadrilaterals,
bounded hyperbolic pyramids and triangular prisms, hyperbolic simplices,
and Lambert cubes in $S^3$, $\E^3$, $\H^3$.

\vspace{5pt}
\noindent
In this paper, we solve this problem for  spherical simplices.

\vspace{20pt}

\noindent
{\bf 2. }
In \cite{simplex}, this problem is solved for  hyperbolic simplices.
The methods of \cite{simplex} are valid for the spherical case, 
and these methods are sufficient to solve the problem for any
partial case. Namely, if $G$ is a reflection group in $S^n$,
then it is possible
to find all the simplices, generating the group $G$.
It is shown in \cite{simplex} that
for any $G$ the procedure takes finite 
(but not uniformly bounded!) time.

A spherical discrete reflection group $G$ is a group generated by 
one of the spherical Coxeter simplices. Any indecomposable spherical 
Coxeter simplex
is one of $A_n$, $B_n$, $D_n$,
$E_6$, $E_7$, $E_8$, $F_4$, $H_3$, $H_4$ and $G_2^{(m)}$
(we use the standard notation, see Table~\ref{coxeter}). 
The rest spherical Coxeter simplices are the direct products of some 
indecomposable Coxeter simplices.
It is sufficient to classify simplices generating indecomposable 
discrete reflection groups (other simplices generating  discrete groups 
are the direct products of these ones).

Since there exist indecomposable Coxeter simplices in any dimension,
it is not possible to solve the general problem case by case.
One could expect that the answer become stable, while $n$ goes to infinity.
But the list of answers grows exponentially.

E.~B.~Vinberg suggested a way to reduce the number of items in 
the answer.
Any simplex in the spherical n-space is bounded by $n+1$ hyperplanes.
Observe, that these hyperplanes decompose the $n$-dimensional sphere
into $2^{n+1}$ simplices. 
 We call this set of simplices  a {\bf family}.
 Evidently, the simplices in a family generate one and the same 
group. We say that a family generates a reflection group.
In this paper we classify all the families of simplices, 
generating discrete groups. 

The classification up to a family is shorter than one for the ordinary
simplices is. The number of families for  $A_n$, $B_n$ and $D_n$
still goes to infinity with growth of the dimension.
But in this situation it is possible to classify the families
in terms of graphs.
In Section 2, we prove that 

\begin{itemize}
\item[$\bullet$]
The families generating $A_n$ 
are in one-to-one correspondence with
the trees with $n+1$ vertices.
 
\item[$\bullet$]
The families generating $B_n$   
are in one-to-one correspondence with
the  trees with $n$ vertices 
exactly one of which is marked.

\item[$\bullet$]
The families generating $D_n$   
are nearby in one-to-one correspondence with
the  connected graphs with $n$ vertices  containing exactly one
cycle (possibly the cycle contains only two vertices, 
but the cycle containing exactly one vertex is permitted).

\end{itemize}

Once the classification for  $A_n$, $B_n$ and $D_n$ is done,
we can solve the problem for the rest groups case by case.
In fact, the answers for some groups are huge even in terms of the families.
For example, there are 78 families generating $H_4$,
223 families generating $E_7$, and 1242 families generating $E_8$.
In Section 3, we derive a criteria in terms of the dihedral angles 
for a family generating a discrete group. 
The complete answer for each of the groups 
$E_6$, $E_7$, $E_8$, $F_4$, $H_3$ and $H_4$ is contained in the
Appendix. (The answer for the group $G_2^{(m)}$
is quite evident: two reflections  generate $G_2^{(m)}$
if and only if the angle between the mirrors equals $\pi \frac{k}{m}$,
where $k$ is prime to $m$).

\subsection*{Acknowledgments}
The work was written during the stay at the University of Fribourg,
Switzerland.
The author would like to thank Prof. R.~Kellerhals
for the invitation and the University
for  the hospitality. 
The author is grateful to Prof. R.~Kellerhals for 
very useful discussions.  
The author is also grateful to Prof. E.~B.~Vinberg
for the idea to classify the simplices modulo family.

\section{Preliminaries}

\noindent
{\bf Spherical reflection groups}.
Let  $S^n$ be the $n$-dimensional sphere, and $G$ be a group acting on $S^n$ 
and generated by 
reflections. Suppose that $G$ acts discretely, that is $G$ is a finite 
group for the spherical case.
Then $G$ is called a Coxeter group and $G$
 is generated by the reflections 
with respect to the facets of a spherical Coxeter polyhedron
(a polyhedron is called Coxeter polyhedron, if its dihedral angles are
the integer parts of $\pi$). 
The classification of the spherical Coxeter polyhedra is due to
Coxeter \cite{Coxeter}.
Any spherical Coxeter polyhedron containing no pair of antipodal points
of $S^n$ is a simplex. Any spherical indecomposable Coxeter simplex
is one of the simplices   $A_n$, $B_n$, $D_n$,
$E_6$, $E_7$, $E_8$, $F_4$, $H_3$, $H_4$ and $G_2^{(m)}$
(a simplex is called indecomposable if it is not a direct product of
some other simplices).

To describe the Coxeter simplices, one can use the Coxeter diagrams.
 The Coxeter diagram of a Coxeter simplex $P$
is a graph whose vertices $v_i$ correspond to the faces $f_i$ of $P$,
the vertices $v_i$ and $v_j$ are joined by a $(k-2)$-fold edge 
if the dihedral angle
formed up by $f_i$ and $f_j$ equals $\frac{\pi}{k}$
(if $f_i$ is orthogonal to $f_j$, $v_i$ and $v_j$ are disjoint).
See Table~\ref{coxeter} for the Coxeter
 diagrams of the indecomposable spherical Coxeter simplices.
For more information about the reflection groups
 see~\cite{29}.

We use one and the same notation for the Coxeter simplex and the group it
generates.

\begin{table}
\begin{center}
\begin{tabular}{cc}
$A_n$ $(n\ge 1)$ & \epsfig{file=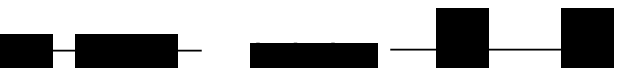,width=0.217\linewidth}
\\
$\vphantom{\int\limits^a}B_n\vphantom{\int\limits^a}$ $(n\ge 2)$ & 
\epsfig{file=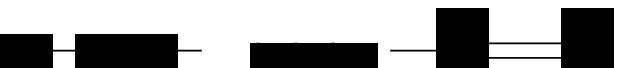,width=0.217\linewidth}
\\
\raisebox{7pt}{$\vphantom{\int\limits^A}D_n\vphantom{\int\limits^A}$ $(n\ge 4)$} & \epsfig{file=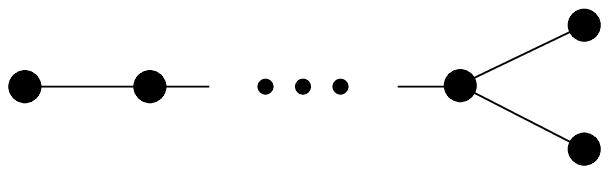,width=0.217\linewidth}
\\
\raisebox{14pt}{$\vphantom{\int^a}E_6\vphantom{\int^a}$} & 
\epsfig{file=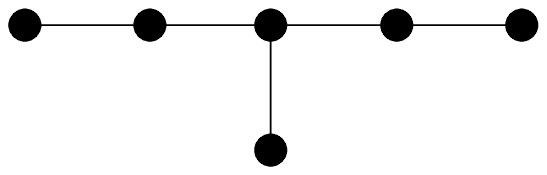,width=0.2\linewidth}
\\
\raisebox{15pt}{$\vphantom{\int^a}E_7\vphantom{\int^a}$} & 
\epsfig{file=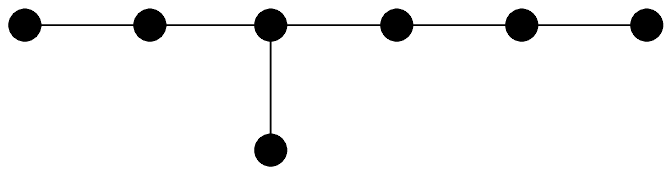,width=0.247\linewidth}
\\
\raisebox{15pt}{$\vphantom{\int^a}E_8\vphantom{\int^a}$} & 
\epsfig{file=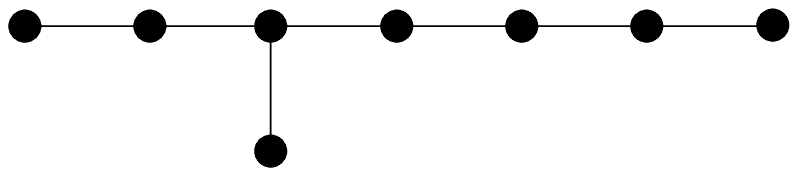,width=0.3\linewidth}
\\
$\vphantom{\int}F_4\vphantom{\int}$ & 
\epsfig{file=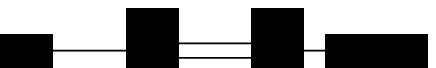,width=0.157\linewidth}
\\
$\vphantom{\int\limits^a}G_2^{(m)}\vphantom{\int\limits^a}$ & 
\psfrag{m}{\scriptsize $m$}
\epsfig{file=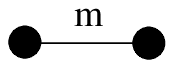,width=0.063\linewidth}
\\
$\vphantom{\int\limits^a}H_3\vphantom{\int\limits^a}$ & 
\epsfig{file=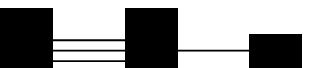,width=0.11\linewidth}
\\
$\vphantom{\int^A}H_4\vphantom{\int^A}$ & 
\epsfig{file=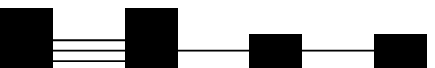,width=0.17\linewidth}
\\
\end{tabular}
\caption{Coxeter diagrams of indecomposable spherical simplices.}
\label{coxeter}
\end{center}
\end{table}

\noindent
{\bf Simplices generating discrete reflection groups}.
Let $\E^{n+1}$ be a $(n+1)$-dimensional Euclidean space and $S^n$
be a unit $n$-sphere centered at the origin. Any hyperplane in $S^n$
corresponds to a hyperplane in  $\E^{n+1}$ containing the origin.
An $n$-dimensional simplex $P$ in $S^n$ corresponds to some $(n+1)$-faced
cone $C$ in  $\E^{n+1}$ with tip at the origin. 
A simplex can be represented by a system of $n+1$ unit vectors 
$f_1,...,f_{n+1}$ orthogonal to the faces of the cone $C$ 
and faced outside of $C$.

Let $\Pi_1,...,\Pi_{n+1}$ be the faces bounding $P$.
The hyperplanes containing $\Pi_1,...,\Pi_{n+1}$
 decompose $S^n$ into $2^{n+1}$ simplices
$P_1$,...,$P_{2^{n+1}}$ encoded by the vectors  
$\{\pm f_1,...,\pm f_{n+1}\}$. 
Throughout this paper the set of simplices $P_1$,...,$P_{2^{n+1}}$
is called a {\bf family}.
Obviously, any of $P_1$,...,$P_{2^{n+1}}$
generates the same reflection group as $P$ does.
Thus, we can study  families instead of ordinary 
simplices. (In fact, any family contains at most $2^n$ simplices 
up to an isometry: 
the simplex $\{f_1,...,f_{n+1}\}$ is always congruent to 
 $\{-f_1,...,-f_{n+1}\}$).

Let $P$ be a simplex generating a discrete reflection group $G$. 
Clearly, the dihedral angles of $P$ are some rational 
numbers multiplied by $\pi$. Moreover, if $G$ is one of
  $A_n$, $D_n$, $E_6$, $E_7$ and  $E_8$, any dihedral angle of $P$
is either right angle or $\frac{\pi}{3}$, or $\frac{2\pi}{3}$.
If $G$ is either $B_n$ or $F_4$, the dihedral angle
can be also $\frac{\pi}{4}$ or $\frac{3\pi}{4}$.  
Observe, that in these cases any dihedral angle is either 
 $\frac{\pi}{k}$ or  $\frac{\pi(k-1)}{k}$,
and that the only indecomposable discrete spherical reflection groups
we omit here are $H_3$ and $H_4$. 

Analogously to the Coxeter simplices, any simplex  
whose dihedral angles are  either 
 $\frac{\pi}{k}$ or  $\frac{\pi(k-1)}{k}$
can be represented by the following {\bf family diagram}:
the vertices $v_i$ correspond to the faces of $P$;
the vertices $v_i$ and $v_j$ are joined by a $(k-2)$-fold edge 
if the angle
formed up by $f_i$ and $f_j$ is either $\frac{\pi}{k}$
or  $\frac{\pi(k-1)}{k}$
(if $f_i$ is orthogonal to $f_j$, $v_i$ and $v_j$ are disjoint).

Clearly, any two simplices in one family
have one and the same family diagram. To specify the simplex in the family,
one can use a {\bf labeled diagram:} we assign an edge with $''-''$
if the correspondent dihedral angle is acute,
and with $''+''$ otherwise. Since the spherical simplices are determined
(up to an isometry)
by their dihedral angles, different simplices have different labeled diagrams.
Note, that some  diagrams and some labeled diagrams
 correspond to no spherical simplex.

We are aimed now to prove that any family diagram corresponds to at most one 
family of spherical simplices with dihedral angles 
 $\frac{\pi}{2}$, $\frac{\pi k}{3}$ or $\frac{\pi l}{4}$.
This means that if $\Sigma$ is a family diagram for a family
$\Phi$
then any labeling of $\Sigma$ either determines  no spherical simplex
or determines a simplex in $\Phi$
(all simplices are considered up to an isometry).

\begin{lemma}
\label{cycle}
Let $\Gamma $ be a cycle, possibly with some 2-fold edges.   
Then $\Gamma$ is a family diagram for
 at most one family of spherical simplices. 

\end{lemma}

\begin{proof}
Let $\pm f_1,...,\pm f_{n+1}$ be a spherical family whose family 
diagram is $\Gamma$.
We can assume that $f_i$ is orthogonal to $f_j$ if $|i-j|>1$ and
$\{i,j\}\ne \{1,n+1\}$.
Fix $f_1$ and choose the sign for $f_2$ to make the angle $\angle f_1 f_2$
non-acute. Then choose the sign for $f_3$ to make the angle $\angle f_2 f_3$
non-acute, and so on. When the sign for $f_{n+1}$ is chosen,
we have two possibilities: either the angle $\angle f_{n+1}f_1$
is acute, or non-acute. Observe, that if  $\angle f_{n+1}f_1$
is non-acute, then  the simplex $P$ determined by  $f_1,...,f_{n+1}$
is acute-angled. Moreover, any dihedral angle of $P$ is either
right angle or  $\frac{\pi}{3}$, or $\frac{\pi}{4}$.
Thus, $P$ is a Coxeter simplex, and the Coxeter
diagram of $P$ cannot be cyclic.
The Coxeter diagram of $P$ coincides with the family diagram of $P$
which is a cycle.
The contradiction shows that  $\angle f_{n+1}f_1$ is acute.

Suppose that  $\Gamma$ corresponds to two different families. 
Choosing the signs as above for both families, we arrive with 
one and the same simplex (since 
the dihedral angles determine spherical
simplex up to an isometry). Thus, the families coincide. 

\end{proof}

\begin{lemma}
\label{unique4}
Let $\Gamma $ be a graph, possibly with some 2-fold edges.
Then $\Gamma$ is a family diagram for
 at most one family of spherical simplices.

\end{lemma}

\begin{proof}
The proof is by induction on the number of vertices of $\Gamma$.
If $\Gamma$ consists of a unique vertex, the lemma is evident.
Assume, inductively, that the lemma is true 
for any graph with at most $n$ vertices.

Let $\Gamma$ be a graph with $n+1$ vertices $v_1,...,v_{n+1}$. 
Suppose that $\Gamma$ corresponds to two different families,
and $f_1,...,f_{n+1}$ and $g_1,...,g_{n+1}$ are the representatives
of these families.
Let $\Gamma^f$ and $\Gamma^g$ be the labeled diagrams for
these simplices. It is sufficient to prove that after some 
changes of signs of vectors $f_1,...,f_{n+1}$
the labeled diagram  $\Gamma^f$ coincides with $\Gamma^g$.
 
Remove from $\Gamma $ the vertex $v_1$ and all the edges ended in $v_1$.
By the assumption, $\Gamma \setminus v_1$ corresponds to at most one 
family. Thus, up to an isometry we have
$f_i=\pm g_i$, $i=2,...,n+1$ (we choose the sign for each $i$ independently). 
Without loss of generality
we can assume that
$f_i=g_i$, $i=2,...,n+1$, and that the labeled diagram 
$\Gamma^f \setminus v_1$ coincides with
$\Gamma^g \setminus v_1$.

Denote the connected components of $\Gamma \setminus v_1$
by $\Gamma_1$,...,$\Gamma_k$. Suppose that $v_1$ is joined with 
$\Gamma_i$ by the edges $e_i^1,...,e_i^{s_i}$.
We can assume that  the edge
$e_1^1$ of $\Gamma^f$  is labeled by ``+''. 
If it is not, we change the sign of every vector
$f_i$ correspondent to some vertex of  $\Gamma_1$.
Note that these changes preserve the labels assigned to the edges 
of $\Gamma_i$, $i=1,...,k$. 
Analogously, we assume that  
$e_1^i$  is labeled by ``+'' for both diagrams $\Gamma^f$ and
 $\Gamma^g$ and for any $i=1,...,k$.
Note, that the labeled subdiagrams $\Gamma^f_i$ 
and $\Gamma^g_i$ still coincide.

Consider the edge $e_1^2$ joining $v_1$ with $\Gamma_1$.
Since $\Gamma_1$ is connected, $e_1^2$ belongs to some cycle in 
$\Gamma_1\bigcup v_1$. By Lemma~\ref{cycle}, any cycle corresponds to at most 
one family of simplices. Thus, the label of $e_1^2$ is determined 
by the labeled diagram $\Gamma_1^f=\Gamma_1^g$.
Therefore, the edge $e_1^2$ has one and the same label in
 $\Gamma^f$ and $\Gamma^g$.
Analogously, the same is true for any edge ended in $v_1$.
Thus, we have  $\Gamma^f=\Gamma^g$, and the lemma is proved. 

\end{proof}

The above lemma is useful for studying of families generating 
 discrete reflection groups distinct from $H_3$ and $H_4$.
To deal with $H_3$ and $H_4$, we extend the definition of the family diagram.
The dihedral angles of any simplex generating $H_3$ or $H_4$ 
belong to the set $\{\frac{\pi}{2},\frac{k\pi}{3},\frac{l\pi}{5}\}$,
where $k=1,2$; $l=1,2,3,4$.
In the {\bf family diagram},
the vertices $v_i$ correspond to the faces of the simplex;
$v_i$ and $v_j$ are joined by a $(k-2)$-fold edge, 
if  $\angle f_if_j$ is either $\frac{\pi}{k}$
or  $\frac{\pi(k-1)}{k}$,
$v_i$ and $v_j$ are joined by a $3$-fold edge
decomposed into 2 parts,
if $\angle f_if_j$ is either $\frac{2\pi}{5}$
or  $\frac{3\pi}{5}$. See Table~\ref{tab-h_4} for some examples.

Clearly, all the simplices in one family have one and the same family diagram.
Unfortunately, the uniqueness of the family correspondent to a diagram 
does not hold now.
For example, the  triangles  
 $(\frac{2\pi}{5}$,$\frac{\pi}{3}$,$\frac{\pi}{3})$
and
 $(\frac{3\pi}{5}$,$\frac{\pi}{3}$,$\frac{\pi}{3})$
have  the same diagram, 
but these triangles belong to different families. 
Nevertheless, usually a family diagram corresponds to at most one family 
generating $H_3$ or $H_4$:

\begin{lemma}
\label{unique5}
Let $\Gamma $ be a family diagram for a family $\Phi$ generating $H_3$ or 
$H_4$.
If $\Gamma $ is not a cycle with 4 vertices, 
then $\Phi$ is the only family 
with family diagram $\Gamma$.
If $\Gamma $ is  a cycle with 4 vertices, 
then $\Gamma$ corresponds to at most two families. 

\end{lemma}

\begin{proof}
First, suppose that $\Phi $ generates $H_3$.
The list of the triangles generating $H_3$ is contained 
in~\cite{Shwarz}. We list the families of the triangles in
Table~\ref{tab-h_3}. It is easy to see that different
families have different family diagrams in this case.

Now, suppose that $\Phi $ generates $H_4$ and $\Gamma$
is not a cycle with 4 vertices.    
In this case, the lemma follows from 
the argument of the proof of Lemma~\ref{unique4} together with the
result of the previous paragraph.

Finally, if $\Gamma $ is a cycle with 4 vertices,
we can assume that three edges correspond to the acute angles of 
the simplex $P\in \Phi$. The rest edge corresponds to either acute
or non-acute dihedral angle. Thus, there are at most two ways 
to specify the dihedral angles of $P$, and $\Gamma$ corresponds to at most two families.

\end{proof}

\noindent
{\bf Remark.}
In Lemma~\ref{unique5}, we do not specify the types of edges in the cycle.
Thus, there could be up to 20 cyclic diagrams corresponding to two
families, generating $H_4$, each. In fact, there are only two such  diagrams.

\section{Simplices generating $A_n$, $B_n$ and $D_n$}

Let $G$ be one of the reflection groups $A_n$, $B_n$ and $D_n$.
We will say that a reflection subgroup $R$ of $G$ is of {\bf maximal rank},
if only the origin is fixed by $R$. 
 Equivalently, $R$ is generated by $n$ reflections,
or $R$ is generated by some $(n-1)$-dimensional simplex.

Let $\Phi$ be a family generating $G$
and $f_1,...,f_n$ be the vectors orthogonal to the faces of $\Phi$. 
Denote by $\Delta$ the correspondent root system $A_n$, $B_n$ or $D_n$.
Then suitably normalized vectors $f_1,...,f_n$ 
belong to $\Delta$.
 Moreover, any linearly independent $n$-tuple of vectors in
 $\Delta$ determines some family generating a maximal rank subgroup of $G$.

To classify all the families generating $G$ it suffices to classify
all the linearly independent $n$-tuples of vectors in $\Delta$
and then eliminate  families generating  proper subgroups of $G$.
Clearly, we are not interested in the ordering of  the vectors 
in the $n$-tuple.
The Weyl group acts on $E^n$ by  reflections, and thus,
it acts on the linearly independent $n$-tuples of vectors in  $\Delta$.
We will not distinguish  $n$-tuples equivalent under this action.

\vspace{8pt}
\noindent
{\bf Definition.}
Let $W$ be a Weyl group of a root system $\Delta$.
We say that a family $\Phi$  is {\bf embedded} into $W$,
if $\Phi$ is represented by some vectors $\pm f_1,...\pm f_n$ chosen from
$\Delta$.

\subsection{Simplices generating $A_n$}

Let $P$ be a simplex generating a group $A_n$.
Then we can assume that
 the vectors $f_1,...,f_n$ belong to 
$\Delta (A_n)=\{\pm(h_i-h_j)\}, 0\le i<j\le n$ (where 
$h_0,...,h_n$ is a standard basis of $\E^{n+1}$).
From the other side, any  linear independent system
 of vectors in $\Delta (A_n)$
corresponds to a simplex generating $A_n$
(according to \cite{Dyn}, $A_n$ has no maximal rank subgroups). 

For any simplex $P=\{f_1,...f_n\}$ generating $A_n$ we construct the following {\bf graph} $\Sigma$:
the vertices $v_0,...,v_n$ of  $\Sigma$ correspond to the vectors 
$h_0,...h_n$, the vertices $v_i$ and $v_j$ are joined by the edge $e_{i,j}$
if one of the vectors $(h_i-h_j)$ and $-(h_i-h_j)$ belong to the set 
$\{f_1,...,f_n\}$. 
Clearly, two simplices in one family have one and the same 
graph. 

To show that the families are in one-to-one correspondence with the graphs,
we need a notion of a dual graph:

\vspace{8pt}
\noindent
{\bf Definition.}
A complete subgraph $C$ of a graph $ \Sigma$
is called {\bf maximal}, if  $ \Sigma$ contains no complete subgraph 
$C_1$ such that $C$ is a subgraph of $C_1$.
  
A graph $ \Sigma^*$ is {\bf dual} to a graph  $ \Sigma$
if the vertices  of  $ \Sigma^*$ correspond to the
maximal complete subgraphs of $\Sigma$, and two vertices
 of $ \Sigma^*$ are joined by a line if and only if they correspond to
 subgraphs in $\Sigma$ having a common vertex.

\vspace{8pt}

If $ \Sigma$ is a tree, the vertices of $ \Sigma^*$ correspond
to the edges of $ \Sigma$, and two vertices are joint 
if the edges are adjacent.

\begin{lemma}
\label{dual}
Let $ \Sigma_1$ and $\Sigma_2$ be the trees with $n+1$ vertices.
If $ \Sigma_1^*=\Sigma_2^*$, then $ \Sigma_1=\Sigma_2$.

\end{lemma}

\begin{figure}[htb]
\begin{center}
\epsfig{file=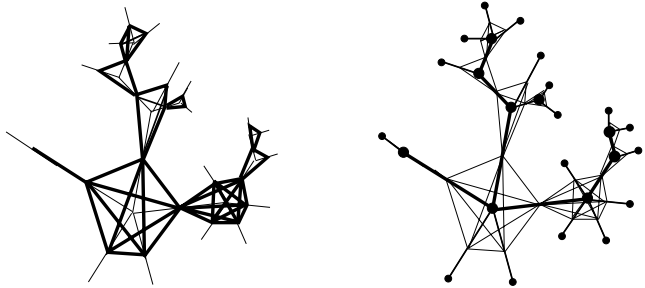,width=0.7\linewidth}
\end{center}
\caption{Cactus and tree.}
\label{fig:cactus}
\end{figure}

\begin{proof}
A graph $\Sigma^*$  dual to a tree looks like a ``cactus'', that is
\begin{itemize}
\item[1.]
 any vertex belongs to at most two maximal complete  subgraphs of $\Sigma^*$;
\vspace{-8pt}
\item[2.]
any pair of maximal complete subgraphs of $\Sigma^*$ has at most one 
common vertex;
\vspace{-8pt}
\item[3.]
any cycle in  $\Sigma^*$ belongs to some complete subgraph.
\end{itemize} 

Any cactus is a graph dual to some tree. 
The lemma follows from the fact that the tree can be reconstructed from 
the cactus by the following algorithm:
\begin{itemize}
\item[1.]
replace each maximal complete subgraph $S$ in $\Sigma^*$ by a vertex $v_S$;
\vspace{-8pt}
\item[2.]
join the vertices $v_{S_1}$ and  $v_{S_2}$, if $S_1$ and $S_2$ have a 
common vertex;
\vspace{-8pt}
\item[3.]
for each vertex $v_S$ attach some additional edges 
to make the valence of $v_S$ equal to the number of vertices in $S$.
\end{itemize} 

See Fig.~\ref{fig:cactus} for an example.   

\end{proof}

\vspace{8pt}
\noindent
{\bf Remark.} 
The proof of Lemma~\ref{dual} shows that $(\Sigma^*)^*$
is  $ \Sigma \setminus L$, where $L$ is the set of all leaves of  
the tree $ \Sigma$ 
(a vertex is called a leaf if its  valency equals 1).

\vspace{8pt}

Lemma~\ref{dual} actually shows that the trees with $n+1$ vertices are
in one-to-one correspondence with the cacti with $n$ vertices.

\begin{theorem}
\label{A}
The families, generating $A_n$, 
are in one-to-one correspondence with the
trees with $n+1$ vertices.

\end{theorem}

\begin{proof}
Let $P=\{f_1,...f_n\}$ be a simplex generating $A_n$, and 
$\Sigma$ be its graph. 
Since $f_0,...,f_n$ are linearly independent, $\Sigma$ has no cycles.
Thus, any family generating $A_n$ corresponds to some tree with 
$n+1$ vertices. 

Let $\Sigma$ be a tree with vertices $v_0,...,v_n$. 
Then $\Sigma$ contains exactly $n$ edges $e_{i,j}$.
Let $f_1,...,f_n$ be the  vectors  
$h_i-h_j, 0\le i<j\le n$, such that $e_{i,j}$ is an edge of $\Sigma$. 
Since  $\Sigma$ is a tree, the vectors $f_1,...,f_n$ are linearly independent.
These vectors define a simplex generating $A_n$.
Thus, the map from the families to the trees is surjective.
It is suffices to prove that this map is injective
and that this map does not depend on the embedding of the family
in $A_n$.

Consider a  graph $\Sigma^*$ dual to $\Sigma$.
Note, that  $\Sigma^*$ coincides with the  family diagram
discussed in the previous section (the vertices correspond to the
vectors $f_1,...f_n$, two vertices are joined if the vectors are 
not mutually orthogonal). 
Since the graphs dual to the trees are in one-to-one correspondence 
with the trees, it is sufficient to prove that
the family diagrams for the families generating $A_n$
are in one-to-one correspondence
with the families generating $A_n$.
The last statement follows from
 Lemma~\ref{unique4}.

\end{proof}


\begin{cor}
Let $\Phi$ be a family of simplices generating $A_n$.
There exists a unique  embedding of $\Phi$ into  $A_n$ (up to the action
of the Weyl group  of $A_n$).

\end{cor}

\begin{proof}
Let $f_1,...f_n$ and $g_1,...,g_n$ be two embeddings of $\Phi$,
and let $\Sigma _f$ and $\Sigma _g$ be the graphs for these embeddings.
Recall, that the vertices of the graphs $\Sigma _f$ and $\Sigma _g$
correspond to the basis vectors $h_0,...,h_n$. 
Theorem~\ref{A} shows that $\Sigma _f$ coincides with  $\Sigma _g$
up to a permutation of vertices.
The Corollary follows from the fact that
the group of permutations of the vectors $h_0,...,h_n$
coincides with the Weyl group of $A_n$. 

\end{proof}

\subsection{Simplices generating $B_n$}

Let $P$ be a simplex generating a group $B_n$.
Then we can assume that the vectors 
 $f_1,...,f_n$ belong to 
$\Delta (B_n)=\{\pm h_i,\ \pm h_i \pm h_j\}, 1\le i<j\le n$ (where 
$h_1,...,h_n$ is a standard basis of $\E^n$).
Any linear independent system of vectors in $\Delta (B_n)$
corresponds to a simplex generating either $B_n$
or some reflection subgroup of $B_n$.

According to~\cite{Dyn}, the only maximal rank indecomposable 
subgroup of $B_n$ is $D_n$.
To generate $B_n$, the system of vectors $f_1,...,f_n$ should be indecomposable
and it should contain at least one of the vectors $\pm h_i$.

For any simplex $P=\{f_1,...f_n\}$ generating $B_n$ we construct the following 
{\bf graph} $\Sigma$:
the vertices $v_1,...,v_n$ of  $\Sigma$ correspond to the vectors 
$h_1,...h_n$,
 the vertices $v_i$ and $v_j$ are joined  
 if $\{f_1,...,f_n\}$ contains either $\pm(h_i+h_j)$ or
 $\pm(h_i-h_j)$
 the vertex $h_i$ is {\bf marked}
 if  $\{f_1,...,f_n\}$ contains $\pm h_i$.
If $\{f_1,...,f_n\}$ contains both $\pm(h_i+h_j)$ and
 $\pm(h_i-h_j)$ then $v_i$ and $v_j$ are joined by two edges. 


Since $f_1,...,f_n$ is indecomposable, the graph $\Sigma$
is connected. Clearly, $\Sigma$ contains at least 
one marked vertex. It is easy to see, that if $\Sigma$ contains two marked
vertices, then the vectors $f_1,...,f_n$ are linearly dependent.
Thus, $\Sigma$ contains a unique marked vertex, and the number of edges in  
 $\Sigma$ is $n-1$.
Since $\Sigma$ has  $n$ vertices,
$\Sigma$ is a tree with one marked vertex
(in particular, no pair of vertices is joined by two edges).
Clearly, any tree with one marked vertex corresponds to a family generating
$B_n$.

To show that the graphs of this type are in one-to-one correspondence
with the families generating $B_n$, we define the {\bf dual} graphs
as the graphs obtained by the following algorithm:
\begin{itemize}
\item[1.]
Construct the graph  $ \Sigma^*$ {\bf dual} to  $\Sigma$
as usually: the vertices of $ \Sigma^*$ correspond
to the edges of $ \Sigma$,  two vertices are joined 
if the edges are adjacent.

\item[2.]
\vspace{-5pt}
Put an additional vertex $v$
in $ \Sigma^*$ for the marked vertex $h_{marked}$ of  $\Sigma$.
\item[3.]
\vspace{-5pt}
Join $v$ by 2-fold line 
with the vertices standing for $\pm h_{marked}\pm h_j$.

\end{itemize}

\begin{lemma}
\label{dualB}
Let $ \Sigma_1$ and $\Sigma_2$ be the trees with $n$ vertices
containing exactly one marked vertex each.
If $ \Sigma_1^*=\Sigma_2^*$ then $ \Sigma_1=\Sigma_2$.

\end{lemma}

\begin{figure}[htb]
\begin{center}
\epsfig{file=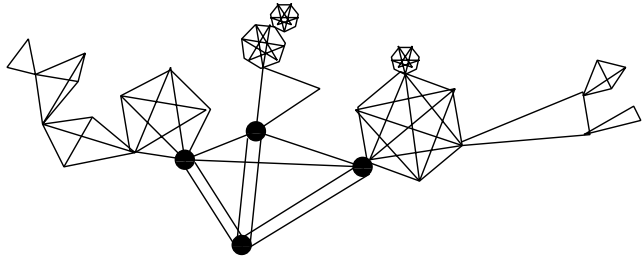,width=0.7\linewidth}
\end{center}
\caption{Bouquet of cacti:
 the additional vertex is joined with
the other vertices by 2-fold lines.
}
\label{fig:bouquet of cacti}
\end{figure}

\begin{proof}
A graph $\Sigma^*$  dual to a tree with one marked vertex 
looks like a ``bouquet  of cacti '',
see Fig.~\ref{fig:bouquet of cacti}.
Any bouquet of cacti is dual to some tree,
and the tree can be easily reconstructed from the bouquet.

\end{proof}

\begin{theorem}
\label{B}
The families generating $B_n$ 
are in one-to-one correspondence with
the  trees with $n$ vertices containing 
 exactly one marked vertex.

\end{theorem}

\begin{proof}
The proof follows  the proof of Theorem~\ref{A}.
It is already shown that
any family generating $B_n$ corresponds to 
a graph described in the theorem, and it is easy to see that any 
graph under consideration corresponds to some family generating $B_n$.
Thus, we have the map from the trees to the families and back.
To see that this map constitutes the one-to-one correspondence,
note, that the dual graph $\Sigma^*$ coincides with the family diagram 
for $f_1,...f_n$, and use Lemma~\ref{dualB} and Lemma~\ref{unique4}.

\end{proof}

\begin{cor}
Let $\Phi$ be a family of simplices generating $B_n$.
There exists a unique embedding of $\Phi$ into  $B_n$ (up to the action
of the Weyl group  of $B_n$).

\end{cor}

\begin{proof}
The proof follows from the fact that
the group of  permutations of the vectors $h_1,...,h_n$
is a subgroup of the Weyl group of $B_n$. 

\end{proof}

\subsection{Simplices generating $D_n$}

Let $P$ be a simplex generating a group $D_n$.
Then we can assume that
the vectors $f_1,...,f_n$ belong to  
$\Delta (D_n)=\{\pm h_i \pm h_j\}, 1\le i<j\le n$ (where 
$h_1,...,h_n$ is a standard basis of $\E^n$).
Any  linear independent system of vectors in $\Delta (D_n)$
corresponds to a simplex generating either $D_n$
or some reflection subgroup of $D_n$.
It follows from~\cite{Dyn}, that any  maximal rank proper reflection subgroup
of $D_n$ is decomposable. Thus, $P$ generates $D_n$
if and only if
 the system of vectors $f_1,...,f_n$ is indecomposable
(the system of vectors  is said to be decomposable if
its Gram matrix is decomposable).

For any simplex $P=\{f_1,...f_n\}$ generating $D_n$ we construct the following 
{\bf colored graph} $\Sigma$:
the vertices $v_1,...,v_n$ of  $\Sigma$ correspond to the vectors 
$h_1,...h_n$, the vertices $v_i$ and $v_j$ are joined by a {\bf black} edge 
 if $\{f_1,...,f_n\}$ contains $\pm(h_i+h_j)$,
$v_i$ and $v_j$ are joined by a {\bf red} edge 
 if $\{f_1,...,f_n\}$ contains $\pm(h_i-h_j)$.
If  $\{f_1,...,f_n\}$ contains both $\pm(h_i+h_j)$ and $\pm(h_i-h_j)$,
$v_i$ and $v_j$ are joined by two edges.

We assume the system $f_1,...,f_n$ to be indecomposable. 
Hence, the graph $\Sigma$
is connected. Since  $\Sigma$ consists of $n$ vertices and $n$ edges,
 $\Sigma$ contains exactly one cycle. The cycle has at least one red
edge, otherwise the vectors  $f_1,...,f_n$ are linearly dependent.
After some changes of signs of $h_1,...,h_n$ one can make all but one
edges of  $\Sigma$ black. Moreover, one can choose any edge of the cycle,
to make this edge the only red edge in $\Sigma$ (the proof follows  the
proof of Lemma~\ref{cycle}). Hence, the graph itself, without colors,
prescribes the family of simplices. The only condition for the uncolored 
connected graph
is that  $\Sigma$ contains a unique cycle (possibly, the cycle contains two
vertices only).

Let $\Sigma$ be a connected graph with exactly one cycle.
Define the {\bf dual} graph $\Sigma^*$ 
as follows:
 the vertices of $ \Sigma^*$ correspond
to the edges of $ \Sigma$,  two vertices are joined 
if the edges have a common vertex;
if two edges make up a cycle,  the correspondent vertices
of $\Sigma^*$ are joined by a dotted edge. 

\begin{lemma}
\label{dualD}
Let $ \Sigma_1$ and $\Sigma_2$ be the graphs with $n$ vertices
containing exactly one cycle each.
If $ \Sigma_1^*=\Sigma_2^*$, then 
 $ \Sigma_1=\Sigma_2$.

\end{lemma}

\begin{proof}
Let $\Sigma$ be a graph with $n$ vertices containing exactly one cycle.
Then $\Sigma^*$  looks like a 
``cactus necklace'',
see Fig.~\ref{fig:cactus necklace}.
If the cycle in $\Sigma$ contains two edges only,
then the string of the necklace is just the dotted edge.
To prove the  lemma, we check that the graph
$\Sigma$ can be recovered from the necklace.

\begin{figure}[htb]
\begin{center}
\psfrag{A}{}
\psfrag{B}{}
\epsfig{file=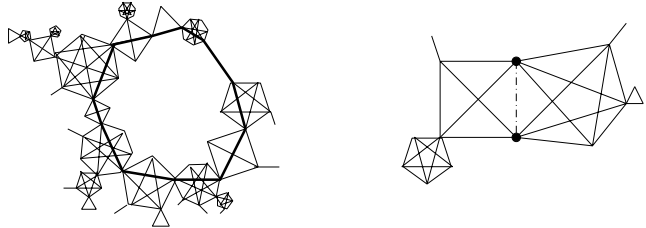,width=0.7\linewidth}
\end{center}
\caption{Cactus necklace. }
\label{fig:cactus necklace}
\end{figure}

Suppose that the string of $\Sigma^*$ contains at least four edges.
Then, to find the graph  $\Sigma$ correspondent to the necklace $ \Sigma^*$,
one can use the same procedure as for $A_n$.

If the string of $\Sigma^*$ consists of two edges,
one can use the usual procedure, but 
the vertices correspondent to two maximal complete subgraphs having 
the dotted edge in common
should be joined by a 2-fold edge.

From now on in this proof
suppose that the string of $\Sigma^*$ consists of three edges.
Then $\Sigma$ could be reconstructed with the usual procedure,
but one should  put no vertex for the string of the necklace
(the string is a complete subgraph with 3 vertices).
Unfortunately, sometimes it is impossible to recognize
the string: see Table~\ref{nondual}, second row, column $ \Sigma_1^*$.
We denote this graph by $T$.

To recognize the string in the other cases, note, that
1)  any edge contained in two different maximal complete subgraphs
belongs to the string;
2) if the string contains at least three vertices and
the maximal complete subgraph $C$ intersects the string,
then $C$ contains exactly two  vertices of the string.
Combining these properties, it is easy to see that
any necklace with at least four vertices contains 
$T$
as a subgraph, and that if the necklace contains at least five
vertices, the string could be recognized.
Thus, the graph $\Sigma$ could be reconstructed,
unless $ \Sigma^*$ is $T$.
The last case does not cause a problem, since $T$ is symmetrical,
and there is no difference which of the two triangles is chosen for a string.
Therefore, $\Sigma$ could be recovered from $ \Sigma^*$,
and the lemma is proved. 

\end{proof}

If $\Sigma^*$ has a dotted edge, denote by $\overline{\Sigma^*}$
the graph $\Sigma^*$ with the dotted edge removed. 
If $\Sigma^*$ contains no  dotted edge, define 
$\overline{\Sigma^*}=\Sigma^*$.

\begin{lemma}
\label{dualD1}
Let $ \Sigma_1$ and $\Sigma_2$ be the graphs with $n$ vertices
containing exactly one cycle each.
If $\overline{ \Sigma_1^*}=\overline{\Sigma_2^*}$, then 
either $ \Sigma_1=\Sigma_2$
or $(\Sigma_1,\Sigma_2)$ is one of two pairs described in 
Table~\ref{nondual}.

\end{lemma}

\begin{table}[htb]
\begin{center}
\begin{tabular}{c}
\psfrag{s1}{$\Sigma_1$}
\psfrag{s2}{$\Sigma_2$}
\psfrag{s1*}{$\Sigma_1^*$}
\psfrag{s2*}{$\Sigma_2^*$}
\psfrag{s*}{$\overline{ \Sigma_1^*}=\overline{\Sigma_2^*}$}
\psfrag{1}{$1$}
\psfrag{2}{$2$}

\epsfig{file=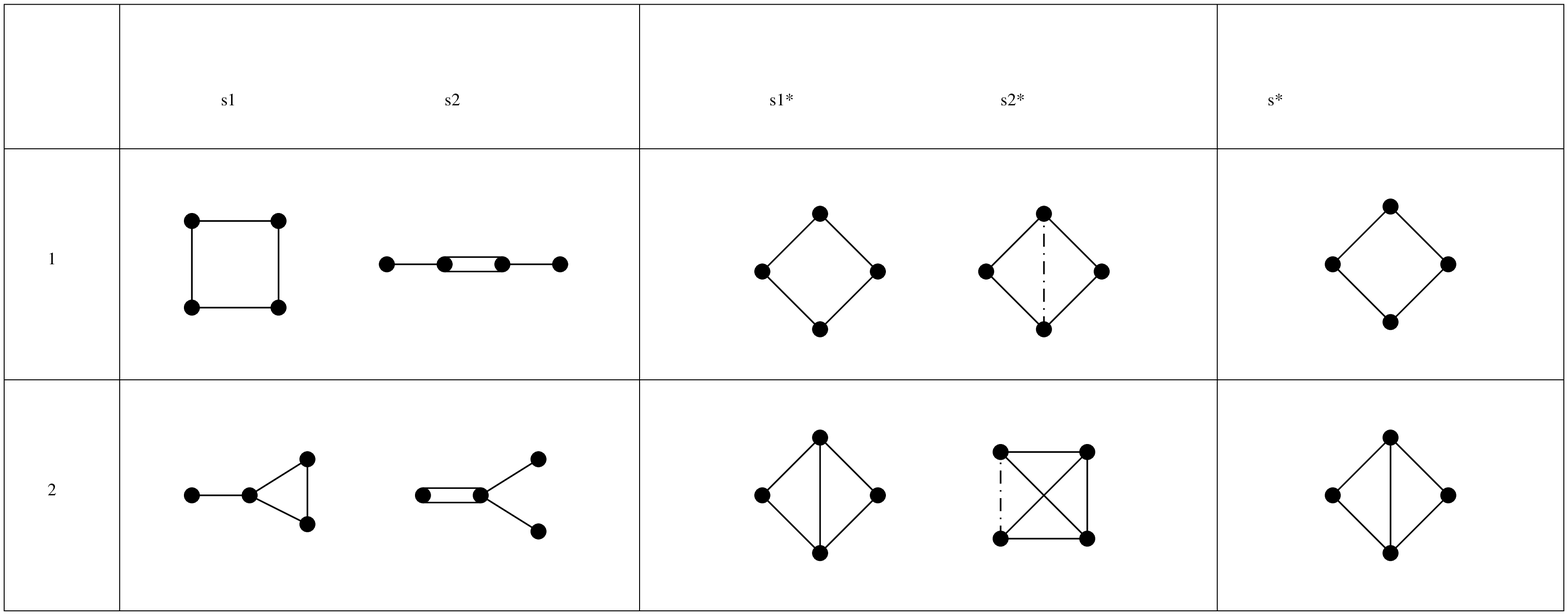,width=0.7\linewidth}\\
\end{tabular}
\end{center}

\caption{Two exclusions.}
\label{nondual}
\end{table}

\begin{proof}
Suppose that neither $ \Sigma_1^*$ nor $\Sigma_2^*$ contain the dotted line.
Then the lemma follows from Lemma~\ref{dualD}.

Suppose that both $ \Sigma_1^*$ and $\Sigma_2^*$ contain the dotted line,
and  $\overline{ \Sigma_1^*}=\overline{\Sigma_2^*}$.
The proof of the preceding lemma shows that the dotted lines in
 $\overline{ \Sigma_1^*}$ and $\overline{\Sigma_2^*}$
could be reconstructed. Thus, if
$\overline{ \Sigma_1^*}=\overline{\Sigma_2^*}$,
then $ \Sigma_1^*=\Sigma_2^*$, and 
the lemma follows from  Lemma~\ref{dualD}.

We are left with the case when
 $ \Sigma_1^*$ contains a dotted line, 
 $\Sigma_2^*$ contains no dotted line,
and  $\overline{ \Sigma_1^*}=\overline{\Sigma_2^*}$.
Suppose that the 2-fold edge in $\Sigma_1$ is incident to 
at least three other edges. Then $\overline{ \Sigma_1^*}$ has at least 
two cycles whose vertices are not the vertices of a complete subgraph.
This is impossible, since  $\overline{ \Sigma_2^*}$ contains 
at most one such a cycle 
(this cycle is a string in $\overline{ \Sigma_2^*}$, 
unless $\Sigma_2=T$, where $T$ is a graph formed up by two triangles).
Suppose that the 2-fold edge in $\Sigma_1$ is incident to 
exactly two other edges. Then $\overline{ \Sigma_1^*}$ has a cycle
$C$ whose vertices are  not the vertices of a complete subgraph.
No edge of $C$ belongs to a complete subgraph having
at least one vertex not in $C$.
This is impossible for $\overline{ \Sigma_2^*}$, unless
any vertex of $\overline{ \Sigma_2^*}$ is contained in $C$.
Thus, $ \Sigma_1$ and  $ \Sigma_2$ contain four vertices each, 
and  $(\Sigma_1,\Sigma_2)$ is one of two pairs described in 
Table~\ref{nondual}.
Finally, suppose that the 2-fold edge in $\Sigma_1$ is incident to 
a unique edge. Then $\overline{ \Sigma_1^*}$ has two leaves with a common
vertex, that is impossible for $\overline{ \Sigma_2^*}$.
(If the 2-fold edge is incident to no other edges, 
the lemma is obvious).

\end{proof}

Denote by $\Gamma_1$ and $\Gamma_2$  the graphs shown in 
the column $ \Sigma_2$
of Table~\ref{nondual}.
  
\begin{theorem}
\label{D}
Let $M$ be the set of the  connected graphs with $n$ vertices 
 containing exactly one
cycle (possibly the cycle contains only two edges, 
but the cycle containing exactly one vertex is permitted).
Then the families, generating $D_n$ 
are in one-to-one correspondence with
the graphs contained in $M\setminus \{\Gamma_1,\Gamma_2 \}$.

\end{theorem}

\begin{proof}
The proof follows the proof of Theorem~\ref{A}.

It is already shown that
any family generating $D_n$ corresponds to some connected 
graph with a unique cycle. 
Any connected graph with $n$ vertices with a unique cycle corresponds to some 
linearly independent indecomposable system of vectors $f_1,...,f_n$
(to construct these vectors take any edge in the cycle, make it red,
and put $h_i-h_j$ for black edges and $h_i+h_j$ for the red one).
Thus, we have the map from the graphs with $n$ vertices and exactly
one cycle 
 to the families generating $D_n$ and back.
To see that this map constitutes the one-to-one correspondence,
note, that the graph $\overline \Sigma^*$ coincides with the family diagram 
for $f_1,...f_n$, and use Lemma~\ref{dualD1} and Lemma~\ref{unique4}.

\end{proof}

Denote by $\Theta_1$ and $\Theta_2$ the graphs shown in the last column of
Table~\ref{nondual}.

\begin{cor}
Let $\Phi$ be a family of simplices generating $D_n$.
If the family diagram for $\Phi$ differs from
$\Theta_1$ and $\Theta_2$, then there exists a unique
 embedding of $\Phi$ into  $D_n$  (up to the action
of the Weyl group  of $D_n$).
If the family diagram for $\Phi$ is either
$\Theta_1$ or $\Theta_2$, then
there are exactly two embeddings.

\end{cor}

\begin{proof}
If the family diagram for $\Phi$ differs from $\Theta_1$ and $\Theta_2$,
the proof follows from the fact that
the group of  permutations of the vectors $h_1,...,h_n$
is a subgroup of the Weyl group of $D_n$. 

Suppose that the family diagram of $\Phi $ is either $\Theta_1$ or $\Theta_2$.
Theorem~\ref{D} shows that $\Phi$ has at most two embeddings
(these embeddings are described by the graphs in the left column of
Table~\ref{nondual}).
To show that these embeddings are different,
note that the Weyl group of $D_4$ never takes the pair 
$\{ h_1-h_2,h_1+h_2\}$ to the pair $\{ h_1-h_2,h_3-h_4\}$.

\end{proof}

\section{Simplices generating other groups}

When all the families generating $A_n$, $B_n$ and $D_n$  
are listed, we are left with finitely many indecomposable spherical reflection
groups. These groups are  
$E_6$, $E_7$, $E_8$, $F_4$, $H_3$ and  $H_4$ 
(we omit $G_2^{(m)}$ since the question is trivial for these groups).
To treat these groups, one can use case-by-case check organized as follows:
\begin{itemize}
\item[1.] 
Let $\Gamma$ be a finite indecomposable reflection group in $S^{n-1}$;
list the unit vectors orthogonal to the planes of the reflections
of the group $\Gamma$ under the consideration (this vectors are collinear to 
the roots in the correspondent root system, unless
$\Gamma_n=H_3$ or  $\Gamma_n=H_4$);
\item[2.] 
for each $n$-tuple of vectors from this list check
that the system is linearly independent and indecomposable;
\item[3.]
each linearly independent $n$-tuple corresponds to some family generating 
$\Gamma$ or some subgroup of $\Gamma$
(but we get every family a lots of times);  
\item[4.]
list the different families obtained at the previous steps
 and eliminate the families generating the proper subgroups of $\Gamma$.
\end{itemize}

Unfortunately, the lists one obtains after all are not too short
(see Appendix for the lists).
To get more compact criteria, we specify the necessary condition
for a family generating a discrete spherical reflection group.

\vspace{5pt}
\noindent
{\bf Notation.} 
Let $P=f_1,...,f_{n+1}$ be a spherical simplex, 
and $\alpha$ be a face of $P$ orthogonal to $f_{n+1}$.
An $(n-1)$-dimensional spherical simplex $f_1,...,f_n$
will be denoted by $P\setminus \alpha$.

\vspace{5pt}
\noindent
{\bf Definition.} 
Let $S$ be a Coxeter simplex generating a reflection group $G$. 
A simplex $P$ {\bf satisfies the subgroup property for a group $G$},
if for any facet $\alpha$ of $P$ there exists a facet $\beta$ of $S$
such that the group generated by $P\setminus \alpha$ 
coincides with the group generated by $S\setminus \beta$.

\begin{lemma}
If a spherical simplex $P$ generates a discrete reflection group $G$,
then $P$ satisfies the subgroup property for $G$. 

\end{lemma}

\begin{proof}
Let $\alpha$ be a face of $P$ and $v$ be a vertex of $P$ opposite to a face
$\alpha$. Then the group generated by $P\setminus \alpha$ 
is a stabilizer $Fix(v)$ of $v$ in $G$.
Let $S$ be a fundamental chamber for $G$ containing $v$,
and $\beta$ be a face of $S$ opposite to $v$. 
Then  $Fix(v)$ is a group generated by
$S\setminus \beta$.

\end{proof}


A case-by-case treating with the computer shows,
that the subgroup property is almost sufficient criteria
for a discreteness of the group generated by $P$: 

\begin{theorem}
\label{subgr}
Let $P$ be a simplex in $S^n$, $n\ge 3$.
Suppose that the  dihedral angles of $P$ are 
in the set $\{\frac{\pi}{2}, \frac{\pi}{3}, \frac{2\pi}{3}, 
\frac{\pi}{4}, \frac{3\pi}{4} \}$.
Let $G_P$ be a reflection group generated by $P$.
Let $G$ be an indecomposable spherical Coxeter group. 

The group $G_P$ is a subgroup of $G$
if and only if $P$  satisfies the subgroup property for the group $G$.

\end{theorem} 

Given a spherical simplex $P$ without the dihedral angles $\frac{k\pi}{5}$,
one can determine whether the  group $G_P$ generated by $P$ is discrete
 or not.
Really, if the dihedral angles of $P$ satisfy the conjecture
of Theorem~\ref{subgr}, one can use the theorem.
Otherwise, $G_P$ is not discrete.  
Note, that if an indecomposable simplex $P$ 
has a dihedral angle $\frac{k\pi}{5}$
and $P$ generates a discrete group $G_P$, then $G_P=H_3$ or $G_P=H_4$.
The analog of Theorem~\ref{subgr} for tetrahedra with a dihedral angle
$\frac{k\pi}{5}$ has some exclusions:

\begin{prop}
\label{5}
Let $P$ be a spherical tetrahedron having a dihedral angle
$\frac{k\pi}{5}$.
Then $P$ generates a discrete reflection group 
if and only if $P$ satisfies the subgroup property for a group
$H_4$ and the family diagram for $P$ differs from the seventeen
diagrams listed in Table~\ref{tab-h_4}.

\end{prop} 

\begin{table}[htb]
\begin{center}
\epsfig{file=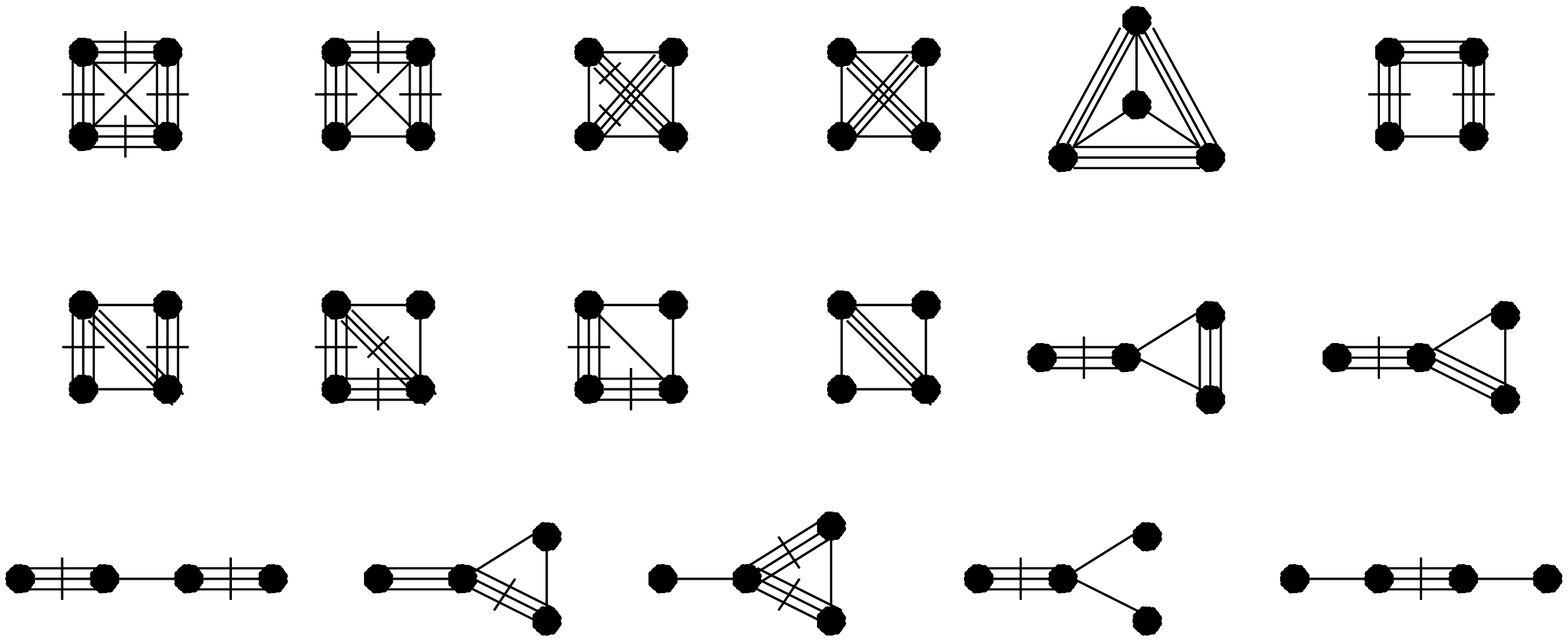,width=0.7\linewidth}
\end{center}
\caption{Family diagrams for spherical tetrahedra, that satisfy the subgroup
property for $H_4$ and generate non-discrete reflection groups.
To represent  dihedral angles $\frac{2\pi}{5}$
and $\frac{3\pi}{5}$ we use 3-fold edge decomposed into 2 parts.  
}
\label{tab-h_4}
\end{table}

\noindent
{\bf Remark.} 
The proof of this proposition as well as of Theorem~\ref{subgr} 
is just a coincidence of the results
of two long case-by-case computer checks.
It is not clear if
 it is possible to prove Theorem~\ref{subgr} a priori.
Why do the exclusions in Proposition~\ref{5} occur,
and how can one find these exclusions a priori?

But one can state that any non-Euclidean simplex whose dihedral angles
belong to the set $\{\frac{\pi}{2}, \frac{\pi}{3}, \frac{2\pi}{3}\}$
generates a discrete reflection group
(see \cite{simplex} for the proof free of case-by-case checks).


\section*{Appendix}
The Appendix contains the list of families generating
$E_6$, $E_7$, $E_8$, $F_4$, $H_3$ and $H_4$.
Families generating $E_6$, $E_7$ and $E_8$ are encoded
in the following way. For a family $\pm f_1,...,\pm f_n$
construct a symmetric matrix $G^+=\{g_{i,j}\}$, where
$g_{i,j}=2|(f_i,f_j)|$.
This is a doubled unsigned Gram matrix of the system
$f_1,..., f_n$.
The upper triangle of $G^+$ is filled up by $0$ and $1$.
Let $p$ be a decimal number which is equal to
the binary number
$\overline{
g_{1,2}g_{1,3}...g_{1,n}g_{2,3}...g_{2,n}...g_{i,i+1}...g_{i,n}...g_{n-1,n}}$.
The number $p$ depends on the ordering of vectors $f_1,...f_n$.
We choose $p$ the smallest possible.
Then $p$ depends only on the family of simplices.
Lemma~\ref{unique4}
shows that different families have different matrices. 
Thus, one can easily recover the family from the number $p$.  
  
\subsection*{Families generating $E_6$}
1187  1195  1197  1199  1214  1259  1260  1261  1263  1270 
1271  1278  1279  1983  3451  3452  3453  3455  3687  3949  
3967  4831  5717  5791  7903  7915  7917  7935  8191  9719  
17598  20287  

\subsection*{Families generating $E_7$}
 
 35050  35051  35052  35053  35054  35055  35070  35071  35098  35099  35102  35103  35195  35196  35197  35199  35431  35693  35711  36208  36209  36211  36212  36213  36215  36282  36283  36286  36287  36518  36519  36653  36663  36724  36725  36735  36862  36863  37117  37461  37523  37535  38547  38551  38623  39510  39595  39599  39647  39661  39666  39679  40670  40683  40685  40694  40703  40895  40956  40959  41191  41463  42487  43503  44967  48127  49279  49342  50303  50367  51515  52031  53055  64511  68986  68987  69487 71245  71247  71262  71265  71295  71339  71403  71675  73451  73577  74063  83070  83259  103864  103865  103867  103871  104012  104013  104015  104052  104053  104055  104063   104103  104107  104171  104178  104179  104371  104383  104444  104445  104447  105038  105039  105059  105071  105321  105339  106345  106349  108015  108295  111081  111085  111437  112527  112545  112557  112575  112639  113643  115838  115839  116027  116028  116029  116031  119341  119423  120510  120511  120637  129023  134739  134893  136925  138975  140765  141015  143327  170577  170589  171739  171740  171741  171743  173533  174555  174559  174799  175071  176013  176095  181437  182815  190103  208715  208719  241484  241487  241495  241502  241535  241546  241579  255574  255575  255607  255663  257774  258047  262143  265551  265558  265559  265711  270055  302991  303095  305639  312487  313719  313775  318975  349535  382463  384511  401101  401119  413021  445695  515577  515807  515820  516095  524287  527486  527487  532095  564863  564927  571646  571647  573375  575550  580159  835063  1050559  1052671  1053407  1083319  1083327  1084269  1084287  1084414  1085438  1086175  1088511  1090039  1150975  1153019  1153021  1160191  1176575   

\subsection*{Families generating $E_8$}
2167311  2167326  2167375  2167378  2167383  2167391  2167415  2167423  2167518  2167531  2167533  2167551  2168396  2168404  2168406  2168414  2168436  2168478  2168492  2168510  2168543  2168556  2168562  2168566  2168575  2168603  2168605  2168607  2168637  2168827  2168829  2168831  2169422  2169578  2170590  2170603  2170605  2170606  2170623  2170733  2170736  2170739  2170748  2170751  2170874  2170879  2171222  2171237  2172399  2172407  2172679  2174823  2175450  2175455  2175465  2175487  2175719  2179199  2179903  2180223  2180287  2180412  2180415  2180927  2182975  2183711  2192383  2193269  2193407  2201234  2201235  2201238  2201239  2201248  2201249  2201251  2201252  2201253  2201255  2201266  2201267  2201270  2201271  2201520  2201521  2201523  2201527  2202212  2202213  2202282  2202283  2202286  2202287  2202348  2202349  2202409  2202413  2202480  2202481  2202491  2202495  2202618  2202619  2202622  2202623  2203570  2203571  2206470  2206471  2206503  2208240  2208241  2208243  2208247  2214202  2214203  2214206  2214207  2216249  2216253  2225013  2226175  2228222  2228223  2231885  2231892  2231893  2231925  2231954  2231955  2231999  2232045  2232093  2234961  2234973  2236123  2236125  2236621  2236631  2236893  2238103  2238943  2267867  2267869  2267871  2267901  2271729  2279611  2279615  2303622  2304906  2305003  2305007  2306015  2306028  2306034  2306047  2308061  2313854  2314043  2314047  2318591  2318934  2320118  2321134  2322431  2323167  2338698  2338700  2338703  2338718  2338721  2338727  2338731  2338733  2338751  2340753  2340765  2352823  2353826  2353827  2353839  2353951  2354926  2355123  2355199  2362710  2362711  2362743  2362799  2362870  2362871  2363143  2363303  2366183  2367207  2370023  2376871  2382335  2402791  2410662  2410663  2415006  2416127  2432931  2433007  2443567  2444662  2444718  2445095  2445551  2445599  2445823  2447787  2479455  2479615  2510061  2510173  2546317  2579967  2580205  2580479  2588671  2612732  2612959  2612972  2613247  2620303  2620415  2621439  2624702  2624703  2624827  2624828  2624829  2624831  2626749  2628127  2639934  2670527  2672703  2694783  2694975  2700735  2702207  2703231  2708031  2709055  2715775  2736062  2890998  2895271  2898151  2932215  2964975  3147775  3148653  3149822  3152895  3180469  3180543  3181567  3182591  3185663  3214188  3214189  3215359  3216108  3216109  3216127  3216379  3218426  3218427  3224575  3228479  3248126  3250171  3251195  3257343  3353583  3367935  3368959  3749823  4161535  4298316  4298317  4298318  4298319  4298336  4298337  4298339  4298347  4298366  4298367  4298410  4298411  4298474  4298475  4298522  4298523  4298537  4298555  4298601  4298611  4298619  4299627  4300648  4300649  4300651  4300655  4302318  4302319  4302695  4305387  4310142  4310143  4310330  4310331  4323181  4329038  4329039  4329059  4329067  4329194  4329195  4329321  4329339  4333291  4335066  4336459  4364890  4368383  4369163  4370251  4370258  4370259  4370273  4370291  4371419  4376127  4376191  4378238  4378427  4382955  4383247  4384335  4384375  4384767  4400619  4408442  4408443  4417103  4417135  4435786  4435791  4435819  4449870  4449891  4449895  4450157  4450175  4452203  4459854  4459855  4462151  4478303  4497231  4497391  4506999  4509935  4576607  4595529  4605211  4607323  4607355  4609369  4642526  4642547  4643563  4649435  4716365  4759167  4767615  4780159  4861307  4908159  5244799  5245803  5249019  5277567  5278571  5280491  5284335  5416927  5419995  6463411  6465420  6465421  6465423  6465463  6465471  6466034  6466035  6466274  6466275  6466379  6466386  6466387  6466401  6466419  6466467  6468307  6470618  6470619  6470655  6473268  6473269  6473271  6473914  6473915  6474040  6474041  6474043  6474047  6474283  6474290  6474291  6478460  6478461  6478463  6478522  6478523  6480319  6480399  6480406  6480407  6480503  6480831  6481487  6481571  6481847  6481919  6482594  6482595  6482803  6488060  6488061  6488063  6496105  6497771  6497775  6498127  6499179  6501193  6502395  6505594  6505595  6505598  6505599  6505787  6506027  6506031  6513019  6513263  6514254  6514255  6514275  6514279  6514541  6515563  6547022  6547023  6547055  6547233  6547237  6547309  6599119  6599141  6605070  6605071  6605091  6607123  6607221  6607277  6607367  6609327  6641255  6673779  6692681  6692685  6702199  6704489  6704496  6704497  6705385  6705741  6706931  6707791  6707947  6739602  6739603  6739699  6740581  6740643  6740655  6740715  6741933  6744839  6749071  6749183  6750207  6782955  6813527  6862783  6863228  6863229  6863661  6869051  6869052  6869053  6869055  6871614  6871615  6877501  6958893  6960819  6960831  6971519  6972983  7005311  7042687  7085671  7091599  7092559  7225175  7274495  7374708  7374709  7374719  7375723  7377777  7377919  7378943  7410537  7410555  7475930  7475931  7475949  7475963  7479785  7479789  7480141  7481343  7482347  7484379  7488127  7497727  7510013  7514103  7517147  7517181  7531518  7531519  8323071  8492625  8492699  8492763  8492785  8496869  8499027  8501213  8559321  8560345  8560349  8562389  8563161  8579807  8580831  8630097  8630109  8644181  8644319  8646381  8670463  8675039  8691405  8691411  8691415  8691423  8694471  8701173  8701591  8702207  8703213  8703325  8703387  8703631  8707807  8724171  8736091  8768748  8768881  8771157  8771231  8781791  8801499  8904415  9027327  9036991  9186783  9472749  9474781 
 9478621$\,$   9478871$\,$  9481183$\,$  9539579$\,$  9539583\quad 10690264\quad 10690265\quad  10690267 \\
 10690271  10692057  10692061  10692301  10692307  10692311  10692319  10693323  10694619  10694623  10699961  10700317  10708565  10708703  10710783  10741395  10741399  10741471  10743518  10743519  10791383  10799265  10799295  10799377  10799515  10800373  10800375  10800383  10803871  10804959  10821323  10821324  10821325  10821327  10821330  10821331  10821342  10821343  10830965  10831027  10831383  10833131  10833138  10833306  10833683  10834139  10834655  10835603  10836557  10836575  10837747  10846175  10865899  10865907  10866611  10868121  10869147  10869343  10869387  10869659  10872727  10873743  10899029  10964223  10964735  10972879  11005654  11005669  11061309  11087291  11091711  11100191  11101375  11124479  11134143  11279957  11280031  11283853  11283935  11284103  11286733  11298903  11570173  11571197  11578335  11602667  11602939  11603967  11604703  11608539  11608543  11609055  11636732  11642869  11873231  12838482  12838483  12838494  12838495  12885710  12885711  12897519  12928107  12928111  12932443  12932683  12932687  12965209  12965453  12965458  12965471  12965473  12965659  12967770  12969542  12969731  12969807  12973907  13028603  13037263  13152127  13380939  13380943  13381455  13418190  13418211  13666029  13666034  13666047  13668347  13669083  13672922  13672923  13673166  13674463  13740011  13740015  15054720  15054723  15054732  15054735  15054750  15054783  15062361  15062364  15062367  15062603  15062605  15062613  15062614  15062623  15062800  15062803  15062812  15062815  15062875  15062877  15063819  15066887  15069953  15069965  15069983  15071093  15134423  15164999  15165007  15165135  15525991  15526127  15526407  15526439  15594739  15594742  15595087  15595175  15595443  15595519  15596111  16033743  16033751  16033783  16043007  16082927  16297919  16303679  16646143  16777215  16881139  16881143  16881223  16881286  16881287  16881319  16881655  16882155  16883535  16883557  16885223  16893206  16893207  16893239  16893303  16893359  16897397  16948951  16991471  17026287  17029711  17095055  17292767  17293031  17308911  17341686  17345959  17348839  17415663  17421431  17422447  17571047  17827319  17860087  19123438  19125911  19126683  19126863  19126926  19126959  19136503  19210215  19214791  19219719  19222759  19222871  19226095  19256663  19257742  19257743  19352974  19357150  19357183  19389919  19438831  19439543  19471599  19471791  19472143  19477903  19478863  19488127  19512823  19518583  19518647  19519599  19519663  19519789  19585183  19611103  19622159  19660287  19668367  19750263  19957239  19958255  19991022  19991399  19993055  19994615  19995631  19997159  20005239  20005295  20531503  21222747  21222991  21353807  21417455  21638479  22054383  22190047  23514531  23517579  23517583  23517662  23517695  23518095  23518686  23518719  23549261  23551467  23779423  23779455  23782679  23789055  24287199  24695071  24695231  25414893  25415255  25415405  25423839  25550303  25611469  25730719  25832653  26318813  26856829  27741421  27741439  28044543  28413951  28876701  28886197  28909278  28909467  28921086  28921087  28921275  28953854  28953855  28954043  28974813  30004351  30141675  30511082  30511083  30983287  30983295  30983422  30983423  33148283  33148286  33148798  33292287  33658042  33658043  33658046  33658047  33658431  33658495  33659002  33659003  33662263  33665143  33665307  33666495  33667007  33670206  33670207  33674815  33761387  33762683  33795871  33795885  33795903  33795967  33809535  33809599  33809975  34000063  34125951  34157695  34638139  34672763  34674303  34707134  34713981  34722367  34777343  35893037  35893055  35893119  35895902  35895927  35895935  35906687  35907127  35907711  35908143  36097215  36175807  36249151  36254847  36254910  36254911  36270143  36388479  36735292  36735295  36739775  36749375  36768059  36769915  36770334  36770335  36771455  36772478  36772667  36779967  36782142  36786751  36805246  36805435  36811133  36841727  36874495  36907387  36921471  37182079  37990187  37992045  37993067  38004271  38189291  38831935  38834815  38966909  38966971  40367358  40367547  40827007  41064124  41068798  41068799  41069375  41070845  41072383  41402367  43093631  45460991  49655295  49657758  49659903  49669567  50572623  50572631  50573063  52944111  52944319  53588471  59236063  66568140  66575871  66584575  67108863  67209076  67209077  67209087  67212011  67212012  67212013  67212283  67213308  67213309  67215855  67275507  67275519  67281373  67281623  67283935  67348447  67348471  67348479  67405303  67438063  67571691  67668799  68191167  68193278  68193279  68289247  68412415  69445599  69445620  69445623  69445631  69636085  69668843  69668862  69668863  69718015  69721087  70253567  70321005  70321023  70325247  70354942  70354943  70355966  70355967  70453247  70454267  70460415  70476799  70517239  70731775  71733227  71733231  71747071  72517627  74358783  74614780  74614783  74615807  74647546  75864029$\,$  76646397$\,$  
87626239$\,$  100902575\quad  100910655\quad  103272255\quad  103792575  \\
134258411  134258413  134258422  134258431  134258623  134258684  134258687  134262695  134265855  134270783  134291179  134291305  134348799  134356703  134361055  134487783  134536703  134749823  134791103  136355563  136355565  136355573  136355574  136355583  136355761  136355767  136355772  136355775  136355836  136355839  136359847  136362999  136363007  136367935  136388461  136394639  136394657  136394669  136394687  136394751  136395755  136402622  136402623  136402749  136411135  136413182  136443903  136444927  136445950  136453851  136453853  136458125  136458207  136472215  136540159  136544255  136577023  136585103  136585207  136666623  136846975  136847039  136862271  137367550  137367551  137370623  137458687  138491883  138540909  138555211  138555215  138604270  138696031  140638079  140638206  140640255  140701618  140705787  140705789  140736507  140771322  140771327  140926700  140926833  140935167  141082495  141082559  141459455  142508031  142745305  142798559  142876365  142876383  142888173  142888285  145055455  149159929  149290975  149290988  149847775  151313911  151575999  153354231  153443823  160964319  168304191$\,$  175521791$\,$  201431037$\,$  204508158$\,$  204508159$\,$  204542971

\subsection*{Families generating $F_4$}

\begin{table}[htb]
\begin{center}
\epsfig{file=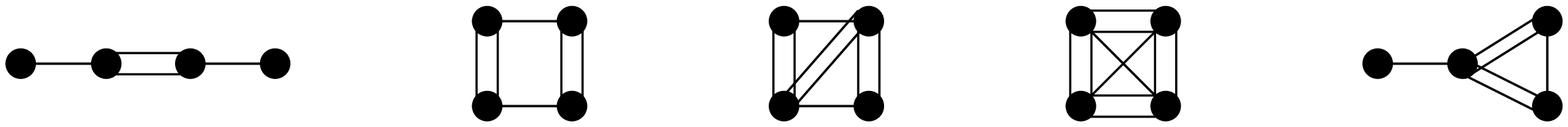,width=0.7\linewidth}
\end{center}
\caption{Families generating $F_4$ are represented by their family diagrams
(in most cases these diagrams are not  Coxeter diagrams).
}
\label{fig:F_4}
\end{table}

\subsection*{Families generating $H_3$}

\begin{table}[htb]
\begin{center}
\begin{tabular}{ccccc}
$(\frac{1}{5},\frac{1}{5},\frac{2}{3})$ &
$(\frac{1}{5},\frac{1}{5},\frac{4}{5})$ &
$(\frac{1}{5},\frac{2}{5},\frac{1}{2})$ &
$(\frac{1}{5},\frac{3}{5},\frac{1}{3})$ &
$(\frac{1}{5},\frac{1}{3},\frac{1}{2})$ \vphantom{$\int_l^A$}\\
$(\frac{1}{5},\frac{1}{3},\frac{2}{3})$ &
$(\frac{2}{5},\frac{2}{5},\frac{2}{5})$ &
$(\frac{2}{5},\frac{1}{3},\frac{3}{5})$ &
$(\frac{2}{5},\frac{1}{3},\frac{1}{3})$ &
$(\frac{2}{5},\frac{1}{3},\frac{1}{2})$ \vphantom{$\int_l^A$}\\

\end{tabular}
\caption{The families of triangles are represented 
by triangles with the smallest
angle sum (and thus, smallest area). The triangle with angles
$\frac{k_1\pi}{l_1},\frac{k_2\pi}{l_2},\frac{k_3\pi}{l_3}$
is denoted by 
$(\frac{k_1}{l_1},\frac{k_2}{l_2},\frac{k_3}{l_3})$.
}
\label{tab-h_3}
\end{center}
\end{table}

\subsection*{Families generating $H_4$}

 Families generating $H_4$ are encoded
in the following way. For a family determined by $\pm f_1,...,\pm f_4$
construct a symmetric matrix $G^+=\{g_{i,j}\}$,
where $g_{i,j}=0$ if $f_i$ is orthogonal to $f_j$,
 $g_{i,j}=1$ if $\angle f_if_j=\frac{\pi}{3}$ or  
$\frac{2\pi}{3}$,
 $g_{i,j}=2$ if $\angle f_if_j=\frac{\pi}{5}$ or  
$\frac{4\pi}{5}$,
and
 $g_{i,j}=3$ if $\angle f_if_j=\frac{2\pi}{5}$ or  
$\frac{3\pi}{5}$.
Let $p$ be a decimal number which is equal to
the base four number
$\overline{
g_{1,2}g_{1,3}...g_{1,n}g_{2,3}...g_{2,n}...g_{i,i+1}...g_{i,n}...g_{n-1,n}}$.
The number $p$ depends on the numbering of vectors $f_1,...f_n$.
We choose $p$ the smallest possible.
Then $p$ depends only on the family of simplices.
Lemma~\ref{unique5}
shows that different families have different matrices,
unless the family diagram is a cycle with 4 vertices. 
Thus, in most cases one can  recover the family from the number $p$.  

\vspace{8pt}
\noindent
  81  85  86  87  91  97  99  101  102  103  106  113  114  117  118  119  127  149  159  163  165  167  177  183  191  213  218  230  234  242  245  246  344  345  346  347  348  349  350  351  364  365  413  414  415  420  422  445  447  489  490  500  503  760  761  1365  1366  1367  1370  1375  1434  1435  1439  1447  1450  1461  1471  1503  1514  1526  1535  1706  1709  1725  1789  2047  2730  4095  

\vspace{8pt}
If the family diagram is a cycle with 4 vertices
and the family generates $H_4$, then $p$ is one of the numbers
344, 348, 364, 420, 500 and 760. It is easy to check that
348 and 500 are the only numbers  
 correspondent to two different families each.

\vspace{25pt}
\noindent
Independent Univ. of Moscow,\\
e-mail: \phantom{ow} felikson@mccme.ru.

\end{document}